\documentclass[a4paper]{article}

\usepackage{amsmath,amsthm}
\usepackage{amssymb,amsfonts}
\usepackage{mathrsfs,color}

\newtheorem{theorem}{Theorem}

\newtheorem{proposition}[theorem]{Proposition}
\newtheorem{lemma}[theorem]{Lemma}
\newtheorem{corollary}[theorem]{Corollary}
\newtheorem{definition}[theorem]{Definition}

\def\N{\mathbb{N}}

\def\R{\mathbb{R}}

\def\<{\langle}
\def\>{\rangle}

\def\ds{\displaystyle} 
\def\div{{\rm div}}
\def\refe#1{(\ref{#1})}

\def\ocirc#1{\ifmmode\setbox0=\hbox{$#1$}\dimen0=\ht0
    \advance\dimen0 by1pt\rlap{\hbox to\wd0{\hss\raise\dimen0
    \hbox{\hskip.2em$\scriptscriptstyle\circ$}\hss}}#1\else
    {\accent"17 #1}\fi}

\def\eps{\varepsilon}
\def\findem{\rule{0.2cm}{0.2cm}}
\def\qed{\findem}

\def\<{\langle}
\def\>{\rangle}
\def\F{\mathcal{F}}
\def\P{\mathbb{P}}
\def\E{\mathbb{E}}
\def\S{\mathbb{S}}
\def\T{\mathbb{T}}
\def\L{\mathscr{L}}
\def\A{\mathcal{A}}
\def\D{D}

\begin{document}

\title{Diffusion limit for a stochastic kinetic problem}
\author{A. Debussche\thanks{IRMAR, ENS Cachan Bretagne, CNRS, UEB. av Robert Schuman, F-35170 Bruz, France. Email: arnaud.debussche@bretagne.ens-cachan.fr} and J. Vovelle\thanks{Universit\'e de Lyon ; CNRS ; Universit\'e Lyon 1, Institut Camille Jordan,  43 boulevard du 11 novembre 1918, F-69622 Villeurbanne Cedex, France. Email: vovelle@math.univ-lyon1.fr}}
\maketitle

\begin{abstract} We study the limit of a kinetic evolution equation involving a small parameter and perturbed by a smooth random term which also involves the small parameter. Generalizing the 
classical method of perturbed test functions, we show the convergence to the solution of a stochastic diffusion equation. 
\end{abstract}

\medskip



\section{Introduction}\label{sec:intro}


Our aim in this work is to develop new tools to study the limit of kinetic equations to fluid models
in the presence of randomness. Without noise, this is a thoroughly studied field in the 
literature. Indeed, kinetic models with small parameters appear in various situations  and it is important to understand the 
limiting equations which are in general much easier to simulate numerically. 

In this article, we consider the following model problem
\begin{equation}
\partial_t f^\eps+\frac{1}{\eps}a(v)\cdot\nabla_x f^\eps=\frac{1}{\eps^2}Lf^\eps+\frac{1}{\eps} f^\eps m^\eps \mbox{ in }\R^+_t\times\T^d_x\times V_v,
\label{spde}\end{equation}
with initial condition
\begin{equation}
f^\eps(0)=f_0^\eps\mbox{ in }\T^d_x\times V_v,
\label{IC}\end{equation}
where $L$ is a {\it linear} operator (see \refe{defL} below) and $m^\eps$ a random process depending on $(t,x)\in\R^+\times\T^d$ (see Section~\ref{s1.1}). We will study the behavior in the limit $\eps\to0$ of its solution $f^\eps$.
\medskip

In the deterministic case $m^\eps=0$, such a problem occurs in various physical situations: we refer to \cite{DegondGoudonPoupaud00} and references therein. The unknown $f^\eps(t,x,v)$ is interpreted as a distribution function of particles, having position $x$ and degrees of freedom $v$ at time $t$. The variable $v$ belongs to a measure space $(V,\mu)$ where $\mu$ is a probability measure. The actual velocity is $a(v)$, where $a\in L^\infty(V;\R^d)$. The operator $L$ expresses the particle interactions. Here, we consider the most basic interaction operator, given by
\begin{equation}
Lf=\int_V f d\mu-f,\quad f\in L^1(V,\mu).
\label{defL}\end{equation}
Note that $L$ is dissipative since
\begin{equation}
-\int_V Lf\cdot f d\mu=\|Lf\|_{L^2(V,\mu)}^2,\quad f\in L^2(V,\mu).
\label{IdL}\end{equation}
In the absence of randomness, the density 
$\rho^\eps=\int_Vf^\eps d\mu$ converges to the solution of the linear parabolic equation (see section \ref{s1.2} 
for a precise statement):
$$
\partial_t\rho-\div(K\nabla\rho)=0\mbox{ in }\R^+_t\times\T^d,
$$
where
\begin{equation}\label{definitionK}
K:=\int_V a(v)\otimes a(v) d\mu(v)
\end{equation}
is assumed to be positive definite. We thus have a diffusion limit in the partial differential equation (PDE) sense.
\medskip

When a random term with the scaling considered here is added to a differential equation, it is classical 
that, at the limit $\eps\to 0$, a stochastic differential equation with time white noise is obtained. This is 
also called a diffusion limit in the probabilistic language, since the solution of such a stochastic
differential equation is generally called a diffusion. Such convergence has been proved initially by Khasminskii \cite{Hasminskii66a, Hasminskii66b} and 
then, using the martingale approach and perturbed test functions, in the 
classical article \cite{PapanicolaouStroockVaradhan77} (see also \cite{EthierKurtz86}, \cite{FouqueGarnierPapanicolaouSolna07}, \cite{Kushner84}).
\medskip

The goal of the present article is twofold. First, we generalize the perturbed test function method to
the context of a PDE and develop some tools for that. We believe that they will be of interest 
for future articles dealing with more complex PDEs. Second, we simultaneously 
take the diffusion limit in the PDE and in the probabilistic sense. This is certainly relevant in a 
situation where a noise with a correlation in time of the same order as a typical length of the
deterministic mechanism is taken. Our main result states that under some assumptions on 
the random term $m$, in particular that it satisfies some mixing properties, the density 
$\rho^\eps=\int_Vf^\eps d\mu$ converges to the solution of the stochastic partial
differential equation
$$
d\rho=\div(K\nabla\rho)dt+ \rho  \circ Q^{1/2} dW(t),\mbox{ in }\R^+_t\times\T^d,
$$
where $K$ is as above, $W$ is a Wiener process in $L^2(\T^d)$ and the covariance 
operator $Q$ can be written in terms of $m$. As is usual in the context of diffusion limit,
the stochastic equation involves a Stratonovitch product.
\medskip

As already mentioned, we use the concept of solution in the martingale sense. This means that 
the distribution of the process satisfies an equation written in terms of the generator (see section
\ref{s2.2} for instance). This generator acts on test functions and the perturbed test function
method is a clever way to choose the test functions such that one can identify the generator
of the limiting equation. Instead of expanding the solution of the random PDE $f^\eps$ as is done
in a Hilbert development in the PDE theory, we work on the test functions acting on the distributions 
of the solutions. 
\medskip

In  section \ref{s2}, we set some notations, describe precisely the random driving term, recall 
the deterministic result and finally state our main result. Section \ref{sec:resolve}
studies
the kinetic equation for $\eps$ fixed. In section \ref{sec:Correctors}, we build the correctors 
involved in the perturbed test function method and identify the limit generator. Finally, in section \ref{sec:DiffLim}, we prove our result. We first show a uniform bound on the $L^2$ norm of the solutions, prove tightness
of the distributions of the solutions and pass to the limit in the martingale formulation.
\medskip

We are not aware of any result on probabilistic diffusion limit using perturbed test functions 
in the context of PDE, but the recent work \cite{deBouardGazeau11} (in a context of nonlinear Schr\"odinger equations) and \cite{PardouxPiatnitski03} (where the underlying PDE is parabolic and the limit $[\eps\to0]$ associated to homogenization effects). A diffusion limit is obtained for the nonlinear Schr\"odinger equation in 
\cite{Marty06}, \cite{DeBouardDebussche10}, \cite{DebusscheTsutsumi10} but there the driving noise is one dimensional and
the solution of the PDE depends continuously on the noise so that in this case an easier argument
can be used. Eventually, note that a method of perturbed test function has also been introduced in the context of viscosity solutions by Evans in \cite{Evans89}. Actually, in the case $m\equiv0$, {\it i.e.} for the deterministic version of \refe{spde}, the method of \cite{Evans89} allows to obtain the diffusive (in the PDE sense) limit $[\eps\to 0]$ of \refe{spde} when the velocity set $V$ is finite.

\section{Preliminary and main result} \label{s2}

\subsection{Notations} 

We work with PDEs on the torus $\T^d$, this means that the space variable $x\in [0,1]^d$ and
periodic boundary conditions are considered. The variable $v$ belongs to a measure space 
$(V,\mu)$ where $\mu$ is a probability measure. We shall write for simplicity  $L^2_{x,v}$
instead of $L^2(\T^d\times V,dx\otimes d\mu)$, its scalar product being denoted by $(\cdot,\cdot)$. We use the same notation for the scalar product of  $L^2(\T^d)$; note that this is consistent since 
$\mu(V)=1$. 
Similarly, we denote by $\|u\|_{L^2}$ the norm $(u,u)^{1/2}$, whether $u\in L^2_{x,v}$ or $L^2(\T^d)$. We use the Sobolev spaces on the torus $H^\gamma(\T^d)$. For $\gamma\in \N$, they consist of periodic functions which are in $L^2(\T^d)$ as well as their derivatives up to order $\gamma$.  For
general $\gamma\ge 0$,  they are easily defined 
by Fourier series for instance. For $\gamma<0$, $H^\gamma(\T^d)$ is the dual of 
$H^{-\gamma}(\T^d)$. 
Classically, for $\gamma_1>\gamma_2$, the injection of $H^{\gamma_1}(\T^d)$ in $H^{\gamma_2}(\T^d)$ is compact.
We use also $L^\infty(\T^d)$ and $W^{1,\infty}(\T^d)$, the subspace of 
$L^\infty(\T^d)$ of functions with derivatives in $L^\infty(\T^d)$. Finally, $L^2(V;H^1(\T^d))$ is the space of functions $f$ of $v$ and $x$ such that all derivatives
with respect to $x$ are in $L^2(\T^d)$ and the square of the norm 
$$
\|f\|_{L^2(H^1)}^2:=\int_V \|f\|_{L^2}^2+\sum_{i=1}^d \|\partial_i f\|_{L^2}^2d\mu
$$
is finite.

\subsection{The driving random term} \label{s1.1}

The random term $m^\eps$ has the scaling 
\begin{equation*}
m^\eps(t,x)=m\left(\frac{t}{\eps^2},x\right),
\end{equation*}
where $m$ is a stationary process  on a probability space $(\Omega,\F,\P)$ and is adapted
 to a filtration $(\F_t)_{t\in\R}$. Note that 
 $m^\eps$ is adapted to the filtration $(\F^\eps_t)_{t\in\R}$, with $\F^\eps_t:=\F_{\eps^{-2}t},\; t\in\R$.
\smallskip

Our basic assumption is that, considered as a random process with values in a space of 
spatially dependent 
functions,  $m$ is a stationary homogeneous Markov process taking values in a subset
$E$ of $W^{1,\infty}(\T^d)$. We assume that $m$ is stochastically continuous. Note that $m$ is supposed not to depend on the variable $v$. The law $\nu$ of $m(t)$ is supposed  to be centered:
\begin{equation}\label{mmoments}
\E m(t)=\int_E n d\nu(n)=0.
\end{equation}
In fact, we also assume that $m$ is uniformly bounded  in $W^{1,\infty}(\T^d)$ so that 
$E$ is included in a ball of $W^{1,\infty}(\T^d)$. We denote by $(P_t)_{t\ge 0}$
a transition semigroup on $E$ associated to $m$ and by $M$ its infinitesimal generator. 
\smallskip

As is usual in the context of diffusion limit, we use the notion of solution of the martingale
problem and need mixing properties on $m$. We 
assume that there is a subset $\D_M$ of $C_b(E)$, the space of bounded continuous functions on $E$, such that, for every $\psi\in\D_M$,
$M\psi$ is well defined  and 
$$
\psi(m(t))-\int_0^t M\psi(m(s))ds 
$$
is a continuous and integrable martingale.
Moreover, we suppose that $m$ is ergodic and satisfies some mixing properties in the sense that there
exists a subspace $\mathscr{P}_M$ of $C_b(E)$ such that  for any $\theta \in  \mathscr{P}_M$ 
the Poisson equation 
\begin{equation}
M\varphi=\theta-\int_E\theta(n) d\nu(n)
\label{InvM}\end{equation}
has a unique solution $\varphi\in \D_M$ satisfying $\ds \int_{E}\varphi d\nu=0$. When $\theta$ satisfies
\begin{equation}
\int_{E}\theta d\nu=0,
\label{conditionInvM}\end{equation}
we denote by $M^{-1}\theta\in\D_M$ this solution and assume that it is given by:
$$
M^{-1}\theta(n)=-\int_0^\infty P_t \theta(n)dt.
$$
In particular, we suppose that the above integral is well defined. It implies that
\begin{equation}\label{mixingPt}
\ds\lim_{t\to+\infty}P_t\theta(n)=0,\quad\forall n\in E.
\end{equation} 
We need that $\mathscr{P}_M$ contains sufficiently many 
functions. In particular, we assume that for each $x\in \T^d$, the evaluation function $\psi_x$ defined by $\psi_x(n)=n(x),\;  n
\in E$,
is in $\mathscr{P}_M$. Also, we assume that, for any $f,g\in L^2_{x,v}$, the function 
$\psi_{f,g}\;:\; n\mapsto (f,ng)$ is in $\mathscr{P}_M$ and we define $M^{-1}I$ from $E$ into $W^{1,\infty}(\T^d)$ by
\begin{equation}\label{defM-1}
(f,M^{-1}I(n)g):=M^{-1} \psi_{f,g}(n),\quad\forall f,g\in L^2_{x,v}.
\end{equation}
We need  that $M^{-1}I$ takes values in a ball of  $W^{1,\infty}(\T^d)$ and take $C_*$
large enough so that
\begin{equation}\label{regNoise}
\|n\|_{W^{1,\infty}(\T^d)}\leq C_*,\quad \|M^{-1}I(n)\|_{W^{1,\infty}(\T^d)}\leq C_*,
\end{equation}
for all $n\in E$. It is natural to require the following compatibility assumption, which would follow from continuity 
properties of $M^{-1}$:
\begin{equation}\label{identityM-1}
M^{-1}\psi_x(n)=M^{-1}I(n)(x),\quad\forall n\in E, x\in\T^d.
\end{equation}
Note that by \eqref{mmoments},  $\psi_{f,g}$ and $\psi_x$ satisfy the centering condition \eqref{conditionInvM}. Note also that, by \refe{defM-1} and \refe{identityM-1}, we have, taking $g=1$,
\begin{equation}\label{identityfM-1}
\int_{\T^d}f(x)M^{-1}\psi_x(n) dx=M^{-1}\psi_{f,1}(n).
\end{equation} 
Eventually, we will also assume that for any $f,g\in L^2_{x,v}$ and for $x\in\T^d$, the functions 
$$
\Psi_{f,g}\;:\; n\mapsto (f,nM^{-1}I(n)g),\quad M^{-1}\psi_{f,1},\quad M^{-1}\psi_x
$$ 
are in $\mathscr{P}_M$.
\medskip

To describe the limit equation, we remark that since $m(0)$ has law $\nu$, 
\begin{align}
-\ds \int_E \psi_y(n)M^{-1}\psi_x(n)d\nu(n)&\ds = -\E\left( \psi_y(m(0))M^{-1}\psi_x(m(0)) \right)\label{expandk0}\\
&\ds = \E\left( \psi_y(m(0))\int_0^\infty P_t\psi_x(m(0))dt \right)\label{ExEy}\\
&\ds = \E\left( \psi_y(m(0))\int_0^\infty \psi_x(m(t))dt \right)\label{Exy}\\
&\ds = \E\left( m(0)(y)\int_0^\infty m(t)(x)dt\right),\label{expandk}
\end{align}
where we have used the Markov property in the identity \refe{ExEy}-\refe{Exy}. We define $k\in L^\infty(\T^d\times\T^d)$ by the formula 
$$
k(x,y)=\E\int_\R m(0)(y)m(t)(x)dt,\quad x,y\in\T^d.
$$
Let $F\in L^\infty(\T^d)$ be the trace
$$
F(x)=k(x,x)=\E\int_\R m(0)(x)m(t)(x)dt,\quad x\in\T^d.
$$
Note that, $m$ being stationary, 
\begin{align}
k(x,y)&\ds = \E\left( \int_0^\infty m(0)(y) m(t)(x)dt\right)+ \E\left( \int_{-\infty}^0 m(0)(y) m(t)(x)dt\right)\nonumber\\
&\ds = \E\left( \int_0^\infty m(0)(y) m(t)(x)dt\right)+ \E\left( \int_{-\infty}^0 m(-t)(y) m(0)(x)dt\right)\nonumber\\
&\ds =  \E\left( m(0)(y)\int_0^\infty m(t)(x)dt\right)+ \E\left( m(0)(x)\int_0^\infty m(t)(y)dt\right),\label{ksymmetric}
\end{align}
so that $k$ is symmetric.  Let $Q$ be the linear operator on $L^2(\T^d)$ associated to the kernel $k$:
$$
Q f(x)=\int_{\T^d} k(x,y) f(y)dy.
$$
\begin{lemma} The operator $Q$ is self-adjoint, compact and non-negative: $(Qf,f)\geq 0$ for all $f\in L^2(\T^d)$.
\label{lem:Qpositive}\end{lemma}

{\bf Proof:} $Q$ is self-adjoint and compact since $k$ is symmetric and bounded. To prove that $(Qf,f)\geq 0$, we will need the following fact: if $\psi\in\mathscr{P}_M$ satisfies \refe{conditionInvM}, then
\begin{equation}\label{PTto0}
\ds\lim_{T\to+\infty}\E|P_T\psi(m(0))|=0.
\end{equation}
Indeed
$$
\E|P_T\psi(m(0))|=\int_E |P_T\psi(n)|d\nu(n)\mbox{ and } |P_T\psi(n)|\leq\|\psi\|_{C_b(E)},
$$
whence \refe{PTto0} by the mixing property \refe{mixingPt} and by the dominated convergence Theorem. In particular, if $\psi\in\mathscr{P}_M$ satisfies \refe{conditionInvM} and if, furthermore, $M^{-1}\psi\in\mathscr{P}_M$, then
$$
P_t\psi=P_tMM^{-1}\psi=\frac{d\;}{dt}P_t M^{-1}\psi,
$$
hence
\begin{equation}\label{PTto00}
\E\left|\int_T^\infty P_t\psi(m(0))\right|=\E|P_TM^{-1}\psi(m(0))|\to 0
\end{equation}
when $T\to+\infty$. For simplicity, let us denote by $\psi_f$ the function $\psi_{f,1}$. By \refe{identityfM-1}, \refe{expandk0}, \refe{ksymmetric} and \refe{PTto00}, we have
\begin{align}
(Qf,f)&=-2\E\left[\psi_f(m(0)) M^{-1}\psi_f(m(0))\right]\nonumber\\
&=2\E\left[\psi_f(m(0)) \int_0^TP_t\psi_f(m(0))dt\right]+o(1),\label{expandQff}
\end{align}
when $T\to+\infty$. On the other hand, for $T>0$, we compute
\begin{align}
\frac{1}{T}\E\left|\int_0^T\psi_f(m(t))dt\right|^2&=\frac{1}{T}\int_0^T\int_0^T \E[\psi_f(m(t)))\psi_f(m(\tau))]dt d\tau \label{expandE21}\\
&=\frac{2}{T}\int_0^T\int_0^t \E[\psi_f(m(t))\psi_f(m(\tau))]dt d\tau\label{expandE22}\\
&=\frac{2}{T}\int_0^T\int_0^t \E[\psi_f(m(t-\tau))\psi_f(m(0))]dt d\tau\label{expandE23}\\
&=\frac{2}{T}\int_0^T\int_0^t \E[\psi_f(m(\tau))\psi_f(m(0))]dt d\tau\nonumber\\
&=\frac{2}{T}\int_0^T(T-\tau) \E[\psi_f(m(\tau))\psi_f(m(0))]d\tau\nonumber\\
&=2\int_0^T\E[\psi_f(m(\tau))\psi_f(m(0))]+r_T\nonumber\\
&=2\E\left[\psi_f(m(0)) \int_0^TP_t\psi_f(m(0))\right]+r_T,
\end{align}
where we have use the homogeneity of $m(t)$ in \refe{expandE22}-\refe{expandE23}. The remainder $r_T$ satisfies
\begin{equation*}
r_T=-2\E\left(\psi_f(m(0))\frac{1}{T}\int_0^T \tau P_\tau\psi_f (m(0)) d\tau\right).
\end{equation*}
Since $\psi_f\in\mathscr{P}_M$, $P_\tau\psi_f=\frac{d\;}{d\tau}P_\tau M^{-1}\psi_f$: this gives
\begin{equation*}
r_T=-2\E\left(\psi_f(m(0))\left[P_TM^{-1}\psi_f(m(0))-\frac{1}{T}\int_0^T P_\tau\psi_f (m(0)) d\tau\right]\right).
\end{equation*}
By \refe{PTto0}, we obtain $r_T=o(1)$. By \refe{expandQff}, $(Qf,f)$ is the limit of the left-hand side of \refe{expandE21}, which is non-negative, hence $(Qf,f)\geq 0$. \qed
\medskip

As a result of Lemma~\ref{lem:Qpositive}, we can define the square root $Q^{1/2}$. 
Note that $Q^{1/2}$ is Hilbert-Schmidt on $L^2(\T^d)$ and that, denoting by $\|Q^{1/2}\|_{\L_2}$ its
Hilbert-Schmidt norm, we have
$$
\|Q^{1/2}\|_{\L_2}^2=\mbox{Tr }Q= \int_{\T^d} k(x,x) dx.
$$
We will not analyze here in detail which kind of processes satisfies our assumptions. The requirement \refe{regNoise} that $m$ and $M^{-1}m$ are a.s. bounded in $W^\infty(\T^d)$ are quite strong. An example of process we may consider is 
$$
m(t)=\sum_{j\in \N} m_j(t)\eta_j
$$
with $\eta_j\in W^{1,\infty}(\T^d)$, 
$$
\sum_{j\in \N} \|\eta_j\|_{W^{1,\infty}(\T^d)} <\infty,
$$
where the processes $(m_j)_{j\in\N}$ are independent real valued 
centered stationary, satisfying the bound 
$$
|m_j(t)|\le C,\; a.s., \; t\in \R,
$$
for a given $C>0$. We are then reduced to analysis on a product space. The invariant measure of $m$ is then easily constructed from the invariant measures of the $m_j$'s. Also, the
Poisson equation can be solved provided each Poisson equation associated to $m_j$ can be 
solved. This can easily be seen by working first on functions $\psi$ depending only on a
finite number of $j$. 

The precise description of the sets $\D_M$ and $\mathscr{P}_M$ depends on the specific processes $m_j,\; j\in\N$.
For instance, if $m_j$ are Poisson processes taking values in finite sets $S_j$, then 
$\D_M$ and $\mathscr{P}_M$ can be taken as the set of bounded functions on $\prod_{i\in\N} S_j$. More general Poisson processes could be considered (see \cite{FouqueGarnierPapanicolaouSolna07}).
\medskip

Actually, the hypothesis~\refe{regNoise} can be slightly relaxed. The boundedness assumption is used two times. First, in the proof of \eqref{e14} and \eqref{e15}, but there it would be sufficient to know that $m$ has finite exponential moments. It is used in a more essential way in Proposition \ref{prop:L2bound}. There, we need that 
the square of the norm of $m$ and $M^{-1}m$ have some exponential moments. However, (under suitable assumptions on the variance of the processes for example), we may consider driving random terms given by Gaussian processes, or more generally diffusion
processes.

\subsection{The deterministic equation}
\label{s1.2}

There are also some structure hypotheses on the first and second moments of $\mu$: we assume
\begin{equation}
\int_V a(v) d\mu(v)=0,
\label{VanishFirstMoment}\end{equation}
and suppose that the following symmetric matrix is definite positive:
\begin{equation}
K:=\int_V a(v)\otimes a(v) d\mu(v) >0.
\label{PositiveSecondMoment}\end{equation}
An example of $(V,\mu,a)$ satisfying the hypotheses above is given by $V=\S^{d-1}$ (the unit sphere of $\R^d$) with $\mu=d-1$-dimensional Hausdorff measure and $a(v)=v$.
\medskip

In the deterministic case $m=0$, the limit problem when $\eps\to0$ is a diffusion equation, as asserted in the following theorem.
\begin{theorem}[Diffusion Limit in the deterministic case] Suppose $m\equiv 0$. Assume that $(f_0^\eps)$ is bounded in $L^2_{x,v}$ and that 
\begin{equation*}
\rho_{0,\eps}:=\int_V f_0^\eps d\mu\to \rho_0\mbox{ in }H^{-1}(\T^d).
\end{equation*}
Assume \refe{VanishFirstMoment}-\refe{PositiveSecondMoment}. Then the density $\rho^\eps:=\int_V f^\eps d\mu$ converges in weak-$L^2_{t,x}$ to the solution $\rho$ to the diffusion equation
\begin{equation*}
\partial_t\rho-\div(K\nabla\rho)=0\mbox{ in }\R^+_t\times\T^d,
\end{equation*}
with initial condition:
$
\rho(0)=\rho_0\mbox{ in }\T^d.
$
\label{th:DiffLimDeterministic}\end{theorem}
This result is a contained is \cite{DegondGoudonPoupaud00} where a more general diffusive limit is analyzed. Note that, actually, strong convergence of $(\rho^\eps)$ can be proved by using compensated compactness, see \cite{DegondGoudonPoupaud00} also.

\subsection{Main result} \label{s1.3}

In our context, the limit of the Problem~\refe{spde}-\refe{IC} is a stochastic diffusion equation. 
\begin{theorem}[Diffusion Limit in the stochastic case] Assume that $(f_0^\eps)$ is bounded in $L^2_{x,v}$ and that 
\begin{equation*}
\rho_{0,\eps}:=\int_V f_0^\eps d\mu\to \rho_0\mbox{ in }L^2(\T^d).
\end{equation*}
Assume \refe{mmoments}-\refe{regNoise}-\refe{VanishFirstMoment}-\refe{PositiveSecondMoment}. 
Then, for all $\eta>0$, the density $\rho^\eps:=\int_V f^\eps d\mu$ converges in law on $C([0,T];H^{-\eta})$ to the solution $\rho$ to the stochastic diffusion equation:
\begin{equation}\label{LimSPDE}
d\rho=\div(K\nabla\rho)dt+\frac12 F\rho + \rho  Q^{1/2} dW(t),\mbox{ in }\R^+_t\times\T^d,
\end{equation}
with initial condition:
$
\rho(0)=\rho_0\mbox{ in }\T^d.
$
In \eqref{LimSPDE}, $W$ is a cylindrical Wiener process on $L^2(\T^d)$.
\label{th:DiffLimStochastic}\end{theorem}

 It is not difficult to see that formally, \refe{LimSPDE} is the It\^o form of 
the  Stratonovitch equation
\begin{equation}\label{LimSPDEStratonovitch}
d\rho=\div(K\nabla\rho)dt+\rho \circ Q^{1/2} dW(t),\mbox{ in }\R^+_t\times\T^d.
\end{equation}

Theorem~\ref{th:DiffLimStochastic} remains true in the slightly more general situation where the coefficient in factor of the noise in \refe{spde} is in the form $\frac{1}{\eps}\sigma(f)m^\eps$ with
\begin{equation*}
\sigma(f)=\bar\sigma(\rho)+f,\quad \rho:=\int_V f d\mu,
\end{equation*}
where $\bar\sigma$ is a smooth, sublinear function.

\section{Resolution of the kinetic Cauchy Problem}\label{sec:resolve}

\subsection{Pathwise solutions}\label{s2.1}

Problem~\refe{spde}-\refe{IC} is linear and solved for instance as follows. Let $A:=a(v)\cdot\nabla_x$ denote the unbounded, skew-adjoint operator on $L^2_{x,v}$ with domain
\begin{equation*}
D(A):=\{f\in L^2_{x,v};a(v)\cdot\nabla_x f\in L^2_{x,v}\}.
\end{equation*} 
Since $A$ is closed and densely defined, by the Hille-Yosida Theorem \cite{CazenaveHaraux98}, it defines a unitary group $e^{tA}$ on $L^2_{x,v}$.

\begin{theorem} Assume \refe{regNoise}. Then, for any $f_0^\eps\in L^2_{x,v}$ and $T>0$, there exists a unique solution $f^\eps$ $\P$-a.s. in $C([0,T];L^2_{x,v})$ of \refe{spde}-\refe{IC} on $[0,T]$, in the sense that, 
\begin{equation*}
f^\eps(t)=e^{-\frac{t}{\eps}A}f_0^\eps+\int_0^t e^{-\frac{t-s}{\eps}A}\left(\frac{1}{\eps^2}Lf^\eps(s)+f^\eps(s)m^\eps(s)\right)ds,
\end{equation*}
$\P$-a.s., for all $t\in [0,T]$. Besides,  if $f_0^\eps\in L^2(V;H^1(\T^d))$, then,  $\P$-a.s. $f^\eps\in C^1([0,T];L^2_{x,v})\cap C([0,T];L^2(V;H^1(\T^d)))$.
\label{th:existssigmabounded}\end{theorem}

The proof of this result is not difficult and left to the reader. The last statement is easily obtained 
since $A$ commutes with derivatives with respect to $x$.
\smallskip

Energy estimates can be obtained. Indeed, for smooth integrable solutions $f^\eps$ to \refe{spde}-\refe{IC}, we have the a priori estimate
\begin{eqnarray*}
\frac{d\;}{dt}\|f^\eps(t)\|_{L^2}^2-\frac{2}{\eps^2}(Lf^\eps,f^\eps)&=&-\frac{2}{\eps}(a(v)\cdot\nabla f^\eps,f^\eps)+\frac{2}{\eps}(f^\eps m^\eps,f^\eps)\\
&=&\frac{2}{\eps}(f^\eps m^\eps,f^\eps).
\end{eqnarray*}
By \refe{IdL} and \refe{regNoise}, this gives the bound
\begin{equation*}
\|f^\eps(t)\|_{L^2}^2+\frac{2}{\eps^2}\int_0^t\|Lf^\eps(s)\|_{L^2}^2 ds\leq \|f_0^\eps\|_{L^2}^2
+\frac{2C_*}{\eps}\int_0^t \|f^\eps(s)\|^2_{L^2_{x,v}} ds,
\end{equation*}
hence, by Gronwall's Lemma, the following bound (depending on $\eps$): 
\begin{equation}
\label{e14}
\|f^\eps(t)\|_{L^2}^2\leq e^{\frac{2C_*}{\eps}t}\|f_0^\eps\|_{L^2}^2.
\end{equation}
Similarly, we have 
\begin{equation}
\label{e15}
\|f^\eps(t)\|_{L^2(H^1)}^2\leq e^{\frac{4C_*}{\eps}t}\|f_0^\eps\|_{L^2(H^1)}^2.
\end{equation}
It is sufficient to assume $f^\eps_0\in L^2(V;H^1(\T^d))$ ({\it resp.} $f^\eps_0\in L^2(V;H^2(\T^d))$) to prove \eqref{e14} ({\it resp.} \eqref{e15}). By density, the inequality holds true for $f^\eps_0\in L^2_{x,v}$ ({\it resp.} $f^\eps_0\in L^2(V;H^1(\T^d))$). In particular, $\|f^\eps(t)\|_{L^2}$ is uniformly bounded in $\omega\in \Omega$ if $f^\eps_0\in L^2_{x,v}$ and $\|f^\eps(t)\|_{L^2(H^1)}$ also if $f^\eps_0\in L^2(V;H^1(\T^d))$.

\subsection{Generator}\label{s2.2}

The process $f^\eps$ is not Markov but the couple $(f^\eps,m^\eps)$ is. 
Its infinitesimal generator is given by : 
\begin{equation}
\L^\eps\varphi=\frac{1}{\eps}\L_{A*}\varphi+\frac{1}{\eps^2}\L_{L*}\varphi,
\label{decLeps}\end{equation}
with
\begin{equation*}\left\{\begin{array}{l }
\L_{A*}\varphi(f,n)=-(Af,D\varphi(f,n))+(f n,D\varphi(f,n)),\\
\\
\L_{L*}\varphi(f,n)=(Lf,D\varphi(f,n))+M\varphi(f,n).
\end{array}\right.\end{equation*}
These are differential operators with respect to the
variables $f\in L^2_{x,v}$, $n\in E$. Here and in the following, $D$ denotes differentiation with respect to $f$ and we identify the differential with the gradient. For a $C^2$ function on $L^2_{x,v}$,
we also use the second 
differential $D^2\varphi$ of a function $\varphi$, it is a bilinear form and we sometimes identify it 
with a bilinear operator on $L^2_{x,v}$, by the formula:
$$
D^2\varphi(f)\cdot(h,k)=(D^2\varphi(f)h,k).
$$
Let us define a set of test functions for the martingale problem associated to the generator $\L^\eps$.
\begin{definition}
We say that $\Psi$ is a good test function if  
\begin{itemize}
\item $\Psi\;:\; L^2(V;H^1(\T^d))\times E\to \R$, $(f,m)\mapsto \Psi(f,m)$ is differentiable with respect 
to $f$  
\item $(f,m)\mapsto D\Psi(f,m)$ is continuous 
from $L^2(V;H^1(\T^d))\times E$ to $L^2_{x,v}$ and maps bounded sets onto bounded sets
\item $(f,m) \mapsto M\Psi(f,m)$ is continuous from 
$L^2(V;H^1(\T^d))\times E$ to $\R$ and  maps bounded sets 
onto bounded sets of $\R$
\item for any $f\in L^2(V;H^1(\T^d))$, $\Psi(f,\cdot)\in \D_M$. 
\end{itemize}
\end{definition}
We have the following result.
\begin{proposition}
\label{t5}
Let $\Psi$ be a good test function. Let $f_0^\eps\in L^2(V;H^1(\T^d))$ and let $f^\eps$ be the solution to Problem~\refe{spde}-\refe{IC}. Then
$$
M^\eps_\Psi(t):=\Psi(f^\eps(t),m^\eps(t))-\int_0^t \L^\eps\Psi(f^\eps(s),m^\eps(s)) ds
$$
is a continuous and integrable $(\F^\eps_t)$ martingale with quadratic variation
\begin{equation}\label{2varM}
\langle M^\eps_\Psi,M^\eps_\Psi\rangle(t)=\int_0^t (\L^\eps|\Psi|^2-2\Psi\L^\eps\Psi)(f^\eps(s),m^\eps(s)) ds.
\end{equation}
\end{proposition}
{\bf Proof:} Let $s,t\ge 0$ and let $s=t_1<\dots<t_n=t$ be a subdivision of $[s,t]$ such that 
$\max_{i}|t_{i+1}-t_i|=\delta$. We have for any $\F^\eps_s$ measurable and bounded $g$
$$
\begin{array}{l}
\E\bigg(\bigg(\Psi(f^\eps(t),m^\eps(t))-\Psi(f^\eps(s),m^\eps(s))\bigg)g\bigg)\\
=\ds \E\bigg( \bigg(\int_s^t  \L^\eps\Psi(f^\eps(\sigma),m^\eps(\sigma)) d\sigma\bigg)g\bigg)+A+B,
\end{array}
$$
With 
\begin{align*}
A\ds &=\sum_{i=1}^{n-1}\E\bigg(\bigg(\Psi(f^\eps(t_{i+1}),m^\eps(t_{i+1}))-\Psi(f^\eps(t_i),m^\eps(t_{i+1}))\\
&\ds -\int_{t_i}^{t_{i+1}}  (-\frac1\eps Af^\eps (\sigma)
+\frac1{\eps^2}Lf^\eps(\sigma)+\frac1\eps f^\eps(\sigma)m^\eps(\sigma),D\Psi(f^\eps(\sigma),m^\eps(\sigma)))d\sigma\bigg)g\bigg)
\end{align*}
and 
$$
\begin{array}{lr}
B&\ds =\sum_{i=1}^{n-1}\E\bigg(\bigg(\Psi(f^\eps(t_{i}),m^\eps(t_{i+1}))-\Psi(f^\eps(t_{i}),m^\eps(t_{i}))\\
&\ds-\int_{t_i}^{t_{i+1}}  M\Psi(f^\eps(\sigma),m^\eps(\sigma))d\sigma\bigg)g\bigg).
\end{array}
$$
We write
\begin{equation*}
A= \E\bigg(\bigg(\int_0^t a_\delta(s)ds\bigg) g\bigg),
\end{equation*}
with
$$
a_\delta(s)\ds = \sum_{i=1}^{n-1}\mathbf{1}_{[t_i,t_{i+1}]}(s) 
\big(  D\Psi(f^\eps(s),m^\eps(t_{i+1}))-D\Psi(f^\eps(s),m^\eps(s))\big)
 \frac{df^\eps}{dt}(s).
$$
Since $f_0^\eps\in L^2(V;H^1(\T^d))$, we deduce from \eqref{e15} and the assumption on 
$\Psi$ that  $a_\delta$ is uniformly integrable with respect to $(s,\omega)$.
Also $f^\eps$ is almost surely continuous and $m^\eps$ is stochastically continuous. It follows that
$\ds D\Psi(f^\eps(s),m^\eps(t_{i+1}))-D\Psi(f^\eps(s),m^\eps(s))$ 
converges to $0$ in probability when $\delta$ goes to zero for any $s$. By uniform integrability, we deduce that 
 $A$ converges to $0$. Similarly, we have
$$
\begin{array}{lr}
B&\ds =\sum_{i=1}^{n-1}\E\bigg(\bigg(\int_{t_i}^{t_{i+1}}  M\Psi(f^\eps(t_i),m^\eps(\sigma))-M\Psi(f^\eps(\sigma),m^\eps(\sigma))d\sigma\bigg)g\bigg),
\end{array}
$$
and, by the same argument, $B$ converges to zero when $\delta$ goes to zero. The result follows~: $M^\eps_\Psi$ is a continuous martingale. Since $\Psi$ is a good test function and $f_0^\eps\in L^2(V;H^1(\T^d))$, it follows from \eqref{e15} and the bound \refe{regNoise} that $t\mapsto\Psi(f^\eps(t),m^\eps(t))$
and $t\mapsto \L^\eps\Psi(f^\eps(t),m^\eps(t))$ are a.s. bounded. The expression~\refe{2varM} for the quadratic variation can then either be computed by expanding
\begin{multline*}
\E|M^\eps_\Psi(t)|^2=\ds\E\left(\left[\sum_{i=1,\dots, n-1} \Psi(f^\eps,m^\eps)(t_{i+1})
- \Psi(f^\eps,m^\eps)(t_i)\right.\right.\\
\left.\left. -\int_{t_i}^{t_{i+1}}\L^\eps\Psi(f^\eps(s),m^\eps(s))ds\right]^2\right),
\end{multline*}
where $0=t_1<\cdots<t_n=t$ is an arbitrary subdivision of $[0,t]$ with step $\delta\downarrow 0$, or, quite similarly, by proceeding
as in Appendix~6.9.1 in \cite{FouqueGarnierPapanicolaouSolna07}. \qed

\section{The limit generator}\label{sec:Correctors}

To prove the convergence of $(\rho^\eps)$, we use the method of the perturbed 
test-function \cite{PapanicolaouStroockVaradhan77}. The method of 
\cite{PapanicolaouStroockVaradhan77} has two steps: first construct a corrector $\varphi^\eps$ to $
\varphi$ so that $\L^\eps\varphi^\eps$ is controlled, then, in a second step, use this with particular test-functions to show the tightness of $(\rho^\eps)$. In the first step, we are led to identify the limit 
generator acting on $\varphi$.

\subsection{Correctors}\label{subsec:Correctors}

In this section, we try to understand the limit equation at $\eps\to 0$. To that purpose, 
we investigate
the limit of the generator $\L^\eps$ by the method of perturbed test-function. 
\smallskip

We restrict our study to smooth test functions and introduce the following class of
functions.
Let $\varphi\in C^3(L^2_{x,v})$. We say that $\varphi$ is regularizing and subquadratic if there exists a constant $C_\varphi\geq 0$ such that  
\begin{equation}\label{regularizingsubquadratic}
\left\{\begin{array}{l}
|\varphi(f)|\leq C_\varphi(1+\|f\|_{L^2})^2,\\ 
\|A^m D\varphi(f)\|_{L^2}\leq C_\varphi(1+\|f\|_{L^2}),\\
|D^2\varphi(f)\cdot (A^{m_1} h,A^{m_2} k)|\leq C_\varphi\|h\|_{L^2}\|k\|_{L^2},\\
|D^3\varphi(f)\cdot (A^{m_1} h,A^{m_2} k,A^{m_3} l)|\leq C_\varphi\|h\|_{L^2}\|k\|_{L^2}\|l\|_{L^2},
\end{array} \right.\end{equation}
for all $f,h,k,l\in L^2_{x,v}$, for all $m,m_i\in\{0,\cdots,3\}$, $i=\{1,2,3\}$.
Note that regularizing and subquadratic functions define good test functions (depending on $f$ only).
\smallskip

Given $\varphi$ regularizing and subquadratic, we want to construct $\varphi_1$, $\varphi_2$ good test functions, such that
\begin{equation*}
\L^\eps\varphi^\eps(f,n)=\L\varphi(f,n)+\mathcal{O}(\eps),\quad \varphi^\eps=\varphi+\eps\varphi_1+\eps^2\varphi_2.
\end{equation*}
The limit generator $\L$ is to be determined. By the decomposition \refe{decLeps}, this is equivalent to the system of equations
\begin{subequations}\label{eq-i}
\begin{align}
\L_{L*}\varphi&=0,\label{eq-2}\\
\L_{A*}\varphi+\L_{L*}\varphi_1&=0,\label{eq-1}\\
\L_{A*}\varphi_1+\L_{L*}\varphi_2&=\L\varphi(f,n),\label{eq0}\\
\L_{A*}\varphi_2&=\mathcal{O}(1).\label{eq1}
\end{align}
\end{subequations}

\subsubsection{Order $\eps^{-2}$}

Equation \refe{eq-2} constrains $\varphi$ to depends on $\rho=\bar f=\int_V f d\mu$ uniquely:
\begin{equation}
\varphi(f)=\varphi(\rho),\quad \rho:=\int_V f d\mu,
\label{varphi}\end{equation}
and imposes that the limit generator $\L$ acts on $\varphi(\bar f)$ uniquely, as expected in the diffusive limit,
in which we obtain an equation on the unknown $\int_V f d\mu$. Indeed, since $\varphi$ is independent on $n$, \refe{eq-2} reads
\begin{equation}
(Lf,D\varphi(f))=0.
\label{eq-21}\end{equation}
Let $(g(t,f))_{t\ge 0}$ denote the flow of $L$ on $L^2(V,\mu)$:
\begin{equation}
\frac{d\;}{dt}g(t,f)=Lg(t,f),\quad g(0,f)=f.
\label{flowL}\end{equation}
An explicit expression for $g$ is 
\begin{equation*}
g(t,f)=\rho+e^{-t}(f-\rho),\quad \rho=\int_V f(v) d\mu(v). 
\end{equation*}
In particular, $g(t,f)\to\rho$ exponentially fast in $L^2(V,\mu)$ when $t\to+\infty$. By \refe{flowL}, equation~\refe{eq-21} is equivalent to 
\begin{equation*}
\varphi(f)=\varphi(g(t,f)),\quad\forall t\in\R,
\end{equation*}
{\it i.e.} \refe{varphi} by letting $t\to+\infty$. 

\subsubsection{Order $\eps^{-1}$}\label{sec:order-1}

Let us now solve the second equation~\refe{eq-1}. To that purpose, we need to invert $\L_{L*}$. Let us work formally in a first step
to derive a solution. Assume that $m(t,n)$ is a Markov process with generator $M$, let $g$ be defined by \refe{flowL} and consider the Markov process $(g(t,f),m(t,n))$. Its generator is precisely 
$\L_{L*}$. Denote by $(Q_t)_{t\ge 0}$ its transition semigroup. Since both $g$ and $m$ satisfy mixing properties, the couple $(g,m)$ also. In particular, we have
\begin{equation}
Q_t\psi(f,n)\to \<\psi\>(\bar f):=\int_{ E} \psi(\bar f,n) d\nu(n),
\label{averageQ}\end{equation}
and it is expected that, under the necessary condition $\<\L_{A*}\varphi\>=0$, a solution to~\refe{eq-1} is given by 
$$
\varphi_1=\int_0^\infty Q_t\L_{A*}\varphi dt.
$$
Let us now compute $\L_{A*}\varphi$. By \refe{varphi}, we have for 
$h\in L^2_{x,v}$, $(h,D\varphi(f))=(\bar h,D\varphi(\rho))$, where as above the upper bar denotes the
average with respect to $v$ and $\rho:=\bar f$. Hence
\begin{equation*}
\L_{A*}\varphi(f,n)=-(\overline{Af},D\varphi(\rho))+(\rho n,D\varphi(\rho)).
\end{equation*}
Since the first moments of $a(v)$ and $m(t)$ vanish, we have
\begin{equation*}
\overline{A\rho}=0\quad\mbox{and}\quad\int_E (\rho n,D\varphi(\rho)) d\nu(n)=0,
\end{equation*}
and the cancellation condition $\<\L_{A*}\varphi\>=0$ is satisfied.  We then write
\begin{align*}
\varphi_1(f,n)&\ds =\int_0^\infty Q_t\L_{A*}\varphi (f,n)dt\\
&\ds=\int_0^\infty \E\left( \L_{A*}\varphi(g(t,f),m(t,n))\right)dt.
\end{align*}
Note that $g$ is deterministic and $\bar g=\rho$, so that
\begin{align*}
\varphi_1(f,n)&\ds = \int_0^{\infty}-(\overline{Ag(t,f)},D\varphi(\rho))+\E\left((\rho m(t,n),D\varphi(\rho))\right)dt\\
&\ds = -\int_0^{\infty}(\overline{Ag(t,f)},D\varphi(\rho))dt-(\rho M^{-1}I(n),D\varphi(\rho)).
\end{align*}
Furthermore, regarding the term $\overline{Ag(t,f)}$, we have
\begin{equation*}
\frac{d\;}{dt}Ag(t,f)=A\frac{d\;}{dt}g(t,f)=ALg(t,f)=A\bar g(t,f)-Ag(t,f).
\end{equation*}
Since $\overline{A\bar f}=0$, we obtain $\ds\frac{d\;}{dt}\overline{Ag(t,f)}=-\overline{Ag(t,f)}$, {\it i.e.} 
$$
\overline{Ag(t,f)}=e^{-t}\overline{Af}.
$$
It follows that
$$
\varphi_1(f,n)=-(\overline{Af},D\varphi(\rho))-(\rho M^{-1}I(n),D\varphi(\rho)).
$$
By \refe{varphi}, this is also equivalent to
\begin{equation}\label{varphi1pre}
\varphi_1(f,n)=-(Af,D\varphi(f))-(f M^{-1}I(n),D\varphi(f)).
\end{equation} 

This computation is formal but it is now easy to define $\varphi_1$ by \eqref{varphi1pre} and
to check that it satisfies~\refe{eq-1}. It is also clear that $\varphi_1$ is a good test function.

\begin{proposition}[First corrector]Let $\varphi\in C^3(L^2_{x,v})$ be regularizing and subquadratic according to \refe{regularizingsubquadratic}. Assume that $\varphi$ satisfy \refe{varphi}. Then \refe{eq-1} has a  solution $\varphi_1\in C^1(L^2_{x,v}\times E)$ given by
\begin{equation}
\varphi_1(f,n)=-(Af,D\varphi(f))-(f M^{-1}I(n),D\varphi(f)),
\label{varphi1}\end{equation}
for all $f\in L^2_{x,v}$, $n\in E$. Moreover $\varphi_1$ is a good test function.
\label{propFirstCorrector}\end{proposition}

\subsubsection{Order $\eps^{0}$}

Let us now analyze Equation \refe{eq0}. Setting $\rho=\bar f$, it gives
\begin{equation}\label{avlimgenerator}
\L\varphi(\rho)=\<\L_{A*}\varphi_1\>(\rho)=\int_E\L_{A*}\varphi_1(\rho,n)d\nu(n).
\end{equation}
We have 
\begin{equation*}
\L_{A*}\psi(f,n)=(-Af+fn,D\psi(f,n))
\end{equation*}
and
\begin{align*}
\varphi_1(f,n)=(-Af-fM^{-1}I(n),D\varphi(f))&=-(\overline{Af},D\varphi(\rho))-(\rho M^{-1}I(n),D\varphi(\rho))\\
&=:\varphi_1^\sharp(f,n)+\varphi_1^*(f,n).
\end{align*}
By \refe{avlimgenerator}, the limit generator is therefore the sum of two terms: 
$$
\L\varphi(\rho)=\L_\sharp\varphi(\rho)+\L_*\varphi(\rho).
$$ 
The first term $\L_\sharp\varphi(\rho)$ corresponds to the deterministic part of the equation. We compute, for $h\in L^2_{x,v}$, 
$$
(h,D\varphi_1^\sharp(f,n))=-(\overline{Ah},D\varphi(\rho))-D^2\varphi(\rho)\cdot(\overline{Af},\bar h).
$$
In particular, evaluating at $f=\rho$ we have
$$
(h,D\varphi_1^\sharp(\rho,n))=-(\overline{Ah},D\varphi(\rho))
$$
since $\overline{A\rho}=0$. Taking then $h=-A\rho+\rho n$ and using once again the cancellation property $\overline{A\rho}=0$, we obtain
$$
\L_\sharp\varphi(\rho)=\int_E (\overline{A^2\rho},D\varphi(\rho))d\nu(n),
$$
{\it i.e.}
\begin{equation}\label{deterministicgenerator}
\L_\sharp\varphi(\rho)=(\overline{A^2\rho},D\varphi(\rho)).
\end{equation}
The second part $\L_*$ corresponds to the random part of the equation: since $\overline{A\rho}=0$,
\begin{equation}\label{randomgenerator}
\L_*\varphi(\rho)=\int_E (\rho n,D\varphi_1^*(\rho,n)) d\nu(n), \quad \varphi_1^*(f,n)=-(\rho M^{-1}I(n),D\varphi(\rho)).
\end{equation}
Now that $\L\varphi=\<\L_{A*}\varphi_1\>$ has been identified, we go on with the resolution of \refe{eq0}. At least formally at a first stage, we can set
$$
\varphi_2(f,n)=-\int_0^\infty Q_t\{\<\L_{A*}\varphi_1\>-\L_{A*}\varphi_1\}(f,n)dt.
$$
To the decomposition $\varphi_1=\varphi_1^\sharp+\varphi_1^*$ then corresponds a similar decomposition 
$$
\varphi_2=\varphi_2^\sharp+\varphi_2^*
$$
for $\varphi_2$. Since $\varphi_1^*(n)$ is linear with respect to $n$, the term
$$
\L_{A*}\varphi_1^*(f,n):=(-Af+fn,D\varphi_1^*(f,n))
$$
can be decomposed into two parts: one that is linear with respect to $n$, the second that is quadratic in $n$. The first (linear) part does not contribute to $\varphi_2^*$ since $m(t)$ is centered: $\E m(t)=0$. Let us thus compute the remaining part
$$
q(f,n):=(fn,D\varphi_1^*(f,n)).
$$
Since $\varphi_1^*(f,n)$ depends on $\rho$ only, we have $q(f,n)=(\rho n,D\varphi_1^*(\rho,n))$. Since, along the flow of $L$, the density $\overline{g(t,f)}=\rho$ is constant, we obtain
$$
\varphi_2^*(f,n)=-\int_0^\infty P_t\left\{\int_E q(\rho,n)d\nu(n)-q(\rho,\cdot)\right\}(n) dt.
$$
In particular, from the expression \refe{randomgenerator} for $\varphi_1^*$ and the fact that $\varphi$ is subquadratic and regularizing, it follows that $\varphi_2^*\in C(L^2_{x,v}\times E)$ is a good test function and satisfies
\begin{align}
|\varphi_2^*(f,n)|\leq C\left(1+\|f\|^2_{L^2}\right),\label{subqvarphi2}\\
|\L_{A*}\varphi_2^*(f,n)|\leq C\left(1+\|f\|^2_{L^2}\right),\label{subqLAvarphi}
\end{align}
for all $f\in L^2_{x,v}$, $n\in E$, where $C$ is a constant depending on the constant $C_*$ in \refe{regNoise} and on the constant $C_\varphi$. Similarly, $\L_{A*}\varphi_1^\sharp(f,n)$ is the sum of one term independent on $n$ and one term linear with respect to $n$. This latter does not contribute to $\varphi_2^\sharp$ by the centering condition $\E m(t)=0$. We explicitly compute the first term:
$$
(-Af,D\varphi_1^\sharp(f,n))=(\overline{A^2 f},D\varphi(\rho))+D^2\varphi(\rho)\cdot(\overline{Af},\overline{Af}).
$$
We have already proved ({\it cf} Section~\ref{sec:order-1}) that, along the flow $g(t,f)$ of $L$, we have $\overline{Ag(t,f)}=e^{-t}\overline{Af}$. Similarly, we have 
$$
\overline{A^2 g(t,f)}-\overline{A^2\rho}=e^{-t}(\overline{A^2 f}-\overline{A^2\rho}).
$$ 
By integrating the exponential $e^{-t}$ with respect to $t$, it follows that
$$
\varphi_2^\sharp(f,n)=(\overline{A^2 f}-\overline{A^2\rho},D\varphi(\rho))+\frac12 D^2\varphi(\rho)(\overline{Af},\overline{Af}).
$$
In particular, $\varphi_2^\sharp\in C(L^2_{x,v}\times E)$ is a good test function and satisfies \refe{subqvarphi2}-\refe{subqLAvarphi}. By introducing
$\A\rho=\div(K\nabla\rho)$, where $K$ is given by \refe{definitionK}, to identify $\L_\sharp$ in \refe{deterministicgenerator}, and by developing the expression \refe{randomgenerator} for $\L_*$, we obtain the following result.
\medskip

\begin{proposition}[Second corrector] Let $\varphi\in C^3(L^2_{x,v})$ be regularizing and subquadratic according to \refe{regularizingsubquadratic}. Assume \refe{varphi} and \refe{mmoments}, \refe{regNoise}, \refe{VanishFirstMoment}, \refe{PositiveSecondMoment}. Let $\A$ denote the unbounded operator defined by
\begin{equation*}
\A\rho=\div(K\nabla\rho),\quad D(\A)=H^2(\T^d)\subset L^2(\T^d).
\end{equation*}
Then \refe{eq0} is satisfied for $\L$ defined by: $\forall\psi\in C^2(L^2(\T^d)),$
\begin{align}
\L\psi(\rho)&=(\A\rho,D\psi(\rho))\nonumber\\
&-\int_E \left\{ (\rho nM^{-1}I(n),D\varphi(\rho))+ D^2\varphi(\rho)\cdot(\rho M^{-1}I(n),\rho n)\right\}d\nu(n),\label{defLlim}
\end{align}
and a corrector $\varphi_2\in C(L^2_{x,v}\times E)$ which is a good test function and satisfies
\begin{align}
|\varphi_2(f,n)|\leq C\left(1+\|f\|^2_{L^2}\right),\nonumber\\
|\L_{A*}\varphi_2(f,n)|\leq C\left(1+\|f\|^2_{L^2}\right),\nonumber
\end{align}
for all $f\in L^2_{x,v}$, $n\in E$, where $C$ is a constant depending on the constant $C_*$ in \refe{regNoise} and on the constant $C_\varphi$.
\label{propSecondCorrector}\end{proposition}
\subsection{Limit equation}\label{sec:limeq}

We will show here that $\L_*$ is the generator of the semi-group associated to a diffusion process on $L^2(\T^d)$. Then \refe{randomgenerator} is a form of $\L_*$ corresponding to the Stratonovitch formulation of the corresponding stochastic differential equation. Actually, we use the expanded form of \refe{defLlim} (which corresponds to a stochastic differential equation in It\^o form) to identify more precisely the limit generator $\L$. The notations for $F$, $k$, $Q$ are those introduced in Section~\ref{s1.1}. We have first:
\begin{align*}
-\int_E  (\rho nM^{-1}I(n),D\varphi(\rho))d\nu(n)
&=\E\int_0^{\infty} (\rho m(0)m(t),D\varphi(\rho))dt\\
&=\frac12\E\int_{\R} (\rho m(0)m(t),D\varphi(\rho))dt\\
&=\frac12 (\rho F,D\varphi(\rho)),
\end{align*}
where 
\begin{equation*}
F(x):=\E\int_{\R}  m(0)(x)m(t)(x)dt=k(x,x).
\end{equation*}
To recognize the part containing $D^2\varphi$, we identify $D^2\varphi$ with its Hessian and
first assume that it is associated to a kernel $\Phi$. Then, we write:
\begin{align*}
&-\int_E D^2\varphi(\rho)\cdot(\rho M^{-1}n,\rho n)d\nu(n)\\
&=\E \int_0^\infty D^2\varphi(\rho)\cdot(\rho m(t),\rho m(0))dt\\
&=\frac12\E \int_{\R} D^2\varphi(\rho)\cdot(\rho m(t),\rho m(0))dt\\
&=\frac12\E \int_{\R} (D^2\varphi(\rho) (\rho m(t)),\rho m(0))dt\\
&= \frac12 \E \int_{\R}\int_{\T^d}\int_{\T^d} \Phi(x,y)\rho(x)m(t)(x)\rho(y)m(0)(y)dx dy dt\\
&=\frac12 \int_{\T^d}\int_{\T^d}\Phi(x,y)k(x,y)\rho(x)\rho(y)dxdy.
\end{align*}
Denote by $q$ the kernel of $Q^{1/2}$, then 
$$
k(x,y)=\int_{\T^d}q(x,z)q(y,z)dz,
$$
which gives
\begin{multline}\label{devLQ}
\int_E D^2\varphi(\rho)\cdot(\rho M^{-1}n,\rho n)d\nu(n)\\
=\frac12 \int_{\T^d}\int_{\T^d}\int_{\T^d}\Phi(x,y)q(x,z)q(y,z)\rho(x)\rho(y)dxdydz\\
=\frac12\mathrm{Trace}[(\rho Q^{1/2})D^2\varphi(\rho)(\rho Q^{1/2})^*].
\end{multline}
By approximation, this formula holds for all $C^2$ function $\varphi$. We conclude that $\L$ is the generator associated to the stochastic PDE
\begin{equation}\label{LimSPDEIto}
d\rho=\div(K\nabla\rho)dt+\frac12 F\rho dt+\rho Q^{1/2}dW(t),
\end{equation}
where $W$ is a cylindrical Wiener process.

\subsection{Summary}\label{sec:sumup}

By Proposition \ref{t5}, Proposition \ref{propFirstCorrector} and Proposition \ref{propSecondCorrector}, we deduce:
\begin{corollary}
\label{c9}
Let $\varphi\in C^3(L^2_{x,v})$ be a regularizing and subquadratic function satisfying \refe{varphi}. There exist two 
good test functions $\varphi_1$, $\varphi_2$ such that, defining $\varphi^\eps= \varphi+\eps\varphi_1+\eps^2\varphi_2$, 
\begin{equation}\label{corrsubq}
\begin{array}{l}
|\varphi_1(f,n)|\leq C\left(1+\|f\|^2_{L^2}\right),\\
\\
|\varphi_2(f,n)|\leq C\left(1+\|f\|^2_{L^2}\right),\\
\\
|\L^\eps\varphi^\eps(f,n)-\L\varphi(f,n)|\leq C\eps\left(1+\|f\|^2_{L^2}\right), 
\end{array}
\end{equation}
for all $f\in L^2_{x,v}$, $n\in E$, where $C$ is a constant depending on the constant $C_*$ in \refe{regNoise} and $C_\varphi$.
Moreover
$$
M^\eps(t):=\varphi^\eps(f^\eps(t),m^\eps(t))-\int_0^t \L^\eps \varphi^\eps(f^\eps(s),m^\eps(s))ds,
\quad t\ge 0,
$$
is a continuous integrable martingale for the filtration
$(\F^\eps_t)$ generated by $m^\eps$ with quadratic variation
\begin{equation}\label{2varM2}
\langle M^\eps,M^\eps\rangle(t)=\int_0^t \left\{(M|\varphi_1|^2-2\varphi_1M\varphi_1)(f^\eps(s),m^\eps(s))+r_\eps(t)\right\}ds,
\end{equation}
where 
\begin{equation}\label{c9remainder}
|r_\eps(t)|\leq C\eps\int_0^t(1+\|f^\eps(t)\|_{L^2}^4)ds,
\end{equation} 
for a constant $C$ depending on $C_*$ and $C_\varphi$. Finally, for $0\le s_1\le \dots \le s_n\le s\le t$ and
$\psi\in C_b((L^2_{x,v})^n)$,
\begin{align}\label{laststatement}
&\ds \left|\E\left(\left(\varphi(\rho^\eps(t))-\varphi(\rho^\eps(s))-\int_s^t\L\varphi(\rho^\eps(\sigma))d\sigma\right)\psi(\rho^\eps(s_1),\dots,\rho^\eps(s_n))\right)\right)\nonumber\\
\nonumber\\
\le&\;\ds C \eps\left(1+\sup_{s\in[0,T]}\E\|f^\eps(t)\|_{L^2}^2\right),
\end{align}
with another constant $C$ depending on the constant $C_*$ in \refe{regNoise}, on $C_\varphi$ and on 
the supremum of $\psi$.
\end{corollary}

{\bf Proof:} Everything has already been proved except for \refe{2varM2} and the last statement~\refe{laststatement}. For this latter, it suffices to 
write:
$$
\begin{array}{l}
\ds\varphi(\rho^\eps(t))-\varphi(\rho^\eps(s))-\int_s^t\L\varphi(\rho^\eps(\sigma))d\sigma\\
\\
\ds =M^\eps(t)-M^\eps(s) -\eps \varphi_1(\rho^\eps(t))-\eps^2\varphi_2(\rho^\eps(t))+\eps \varphi_1(\rho^\eps(s))
+\eps^2\varphi_2(\rho^\eps(s))\\
\\
\ds -\int_s^t\bigg(\L\varphi(\rho^\eps(\sigma))-\L^\eps\varphi^\eps(\rho^\eps(\sigma))\bigg)
d\sigma.
\end{array}
$$
Then, we multiply by $\psi(\rho^\eps(s_1),\dots,\rho^\eps(s_n))$, take the expectation and use the bounds \refe{corrsubq} to conclude. Furthermore, 
\begin{equation}\label{transparent2var}
\mathcal{M}|\varphi|^2-2\varphi \mathcal{M}\varphi=0
\end{equation}
if $\varphi\mapsto\mathcal{M}\varphi$ is a linear first order operator in $\varphi$. Applying \refe{transparent2var} to 
$$
\mathcal{M}\varphi(f,n)=\frac{1}{\eps}\L_{A*}\varphi(f,n)+\frac{1}{\eps^2}(Lf,D\varphi(f,n))
$$
gives 
$$
\L^\eps|\varphi^\eps|^2-2\varphi^\eps\L^\eps\varphi^\eps=M|\varphi_1|^2-2\varphi_1M\varphi_1+r_\eps.
$$ 
By \refe{corrsubq}, the remainder $r_\eps$ satisfies \refe{c9remainder}. \qed

\section{Diffusive limit}\label{sec:DiffLim}

Our aim now is to prove the convergence in law of $\rho^\eps=\int_V f^\eps d\mu$ to $\rho$, solution to  
\refe{LimSPDEStratonovitch}, or equivalently  Equation~\refe{LimSPDEIto}. To that purpose, we
use again the perturbed test function method to get a bound 
on the solutions in $L^2_{x,v}$ then we prove that $\rho^\eps$ is tight in $C([0,T];H^{-\eta})$, $\eta>0$. We follow (and adapt to our context) the method in \cite{FouqueGarnierPapanicolaouSolna07}, paragraph~6.3.5. In particular, we use Kolmogorov criterion to get tightness in section~\ref{sec:tight}; an alternative method would be to use Aldous' criterion for tightness ({\it e.g.} Theorem~4.5 in \cite{JacodShiryaev03}).

\subsection{Bound in $L^2_{x,v}$}\label{sec:L2bound}

\begin{proposition}[Uniform $L^2_{x,v}$ bound] Assume \refe{regNoise}. Then, for all $T>0$, $p\geq 1$, we have
$$
\E\sup_{t\in[0,T]}\|f^\eps(t)\|_{L^2}^p\leq C
$$
where the constant $C\geq 0$ depends on $T$, on $p$, on $\|a\|_{L^\infty(V)}$, on the constant $C_*$ in \refe{regNoise} and $\sup_{\eps>0}\|f^\eps_0\|_{L^2}$ only.
\label{prop:L2bound}\end{proposition}
\medskip

{\bf Proof:} Fix $p\geq 2$. Let us write $a(\eps,t)\lesssim b(\eps,t)$ if there exists a constant $C$ depending on $T$, on $p$, on $\|a\|_{L^\infty(V)}$ and on the constant $C_*$ in \refe{regNoise} only such that $a(\eps,t)\leq C b(\eps,t)$ for all $t\in[0,T]$. Set $\varphi(f):=\frac{1}{2}\|f\|_{L^2}^2$. We want to apply Corollary~\ref{c9} to $\varphi$. This requires some care since $\varphi$ is actually a function of $f$ and not of $\rho$. Thus, we first seek for one corrector $\varphi_1\in C^2(L^2_{x,v}\times E)$ such that, for the modified test-function 
\begin{equation*}
\varphi^\eps:=\varphi+\eps\varphi_1,
\end{equation*}
the term $\L^\eps\varphi^\eps(f^\eps,m^\eps)$ can be accurately controlled: for $f\in L^2(V;H^1(\T^d))$, $n\in E$, we compute
\begin{equation}\label{expandLeps}
\L^\eps\varphi^\eps(f,n)=\eps^{-2}\L_{L*}\varphi(f)+\eps^{-1}(\L_{A*}\varphi+\L_{L*}\varphi_1)(f,n)+\L_{A*}\varphi_1(f,n).
\end{equation}
Since $M\varphi(f,n)=0$ ($\varphi$ being independent on $n$), and since $D\varphi(f,n)=f$, the first term in \refe{expandLeps} is
\begin{equation}
\eps^{-2}\L_{L*}\varphi(f)=-\frac{1}{\eps^2}\|Lf\|_{L^2}^2,
\label{boundI0}\end{equation}
which has a favorable sign. Since $A$ is skew-symmetric, $\L_{A*}\varphi(f,n)=(fn,f)$. This term is difficult to control and we choose $\varphi_1$ to compensate it. We set 
\begin{equation*}
\varphi_1(f,n)=-(fM^{-1}I(n),f),
\end{equation*} 
so that $M\varphi_1=-(fn,f)$ and
the second term in \refe{expandLeps} is
\begin{equation*}
\eps^{-1}(\L_{A*}\varphi+\L_{L*}\varphi_1)=\frac{1}{\eps}(Lf,D\varphi_1(f,n))=-\frac{2}{\eps}(Lf,fM^{-1}I(n)).
\end{equation*}
By \refe{regNoise}, we obtain
\begin{equation}
\eps^{-1}(\L_{A*}\varphi+\L_{L*}\varphi_1)(f^\eps(t),m^\eps(t))\leq \frac{1}{4\eps^2}\|Lf^\eps(t)\|_{L^2}^2+4C_*^2\|f^\eps(t)\|_{L^2}^2.
\label{boundI1}\end{equation}
The remainder $\L_{A*}\varphi_1$ satisfies the following bounds
\begin{align*}
\L_{A*}\varphi_1(f,n)&=-(Af,fM^{-1}I(n))+(fn,fM^{-1}I(n))\\
&=\frac12(f^2,AM^{-1}I(n))+(fn,fM^{-1}I(n))\\
&\leq \left(\frac12\|AM^{-1}I(n)\|_{L^\infty_{x,v}}+\|M^{-1}I(n)\|_{L^\infty_{x,v}}\|n\|_{L^\infty_{x,v}}\right)\|f\|_{L^2}^2\\
&\leq \left(\frac12\|a\|_{L^\infty_v}\|M^{-1}I(n)\|_{W^{1,\infty}_{x}}+\|M^{-1}I(n)\|_{L^\infty_{x,v}}\|n\|_{L^\infty_{x,v}}\right)\|f\|_{L^2}^2.
\end{align*}
By \refe{regNoise}, \refe{boundI0}, \refe{boundI1}, we obtain
\begin{equation}\label{boundlepsabove}
\L^\eps\varphi^\eps(f^\eps(t),m^\eps(t))\lesssim\|f^\eps(t)\|_{L^2}^2.
\end{equation}
Set 
$$
M^\eps(t):=\varphi^\eps(f^\eps(t),m^\eps(t))-\varphi^\eps(f^\eps_0,m^\eps(0))-\int_0^t \L^\eps \varphi^\eps(f^\eps(s),m^\eps(s))ds.
$$
By definition of $\varphi$, $\varphi^\eps$ and $M^\eps$, we have
\begin{multline*}
\frac12\|f^\eps(t)\|_{L^2_{x,v}}^2
=\frac12\|f^\eps_0\|_{L^2_{x,v}}^2-\eps(\varphi_1(f^\eps(t),m^\eps(t)))-\varphi_1(f^\eps_0,m^\eps(0)))\\
+\int_0^t\L^\eps\varphi^\eps(f^\eps(s),m^\eps(s))ds+M^\eps(t).
\end{multline*}
By \refe{boundlepsabove} and the estimate 
\begin{equation}\label{quadvarphi1l2}
|\varphi_1(f^\eps(t),m^\eps(t))|\lesssim \|f^\eps(t)\|_{L^2_{x,v}}^2, 
\end{equation}
we deduce the bound
$$
\|f^\eps(t)\|_{L^2_{x,v}}^2\lesssim \|f^\eps_0\|_{L^2_{x,v}}^2+\eps\|f^\eps(t)\|_{L^2_{x,v}}^2+\int_0^t \|f^\eps(s)\|_{L^2_{x,v}}^2ds+\sup_{t\in[0,T]}|M^\eps(t)|.
$$
For $\eps$ small enough, it follows that
$$
\|f^\eps(t)\|_{L^2_{x,v}}^2\lesssim \|f^\eps_0\|_{L^2_{x,v}}^2+\int_0^t \|f^\eps(s)\|_{L^2_{x,v}}^2ds+\sup_{t\in[0,T]}|M^\eps(t)|,
$$
and, by Gronwall Lemma,
\begin{equation}\label{boundl2Meps}
\|f^\eps(t)\|_{L^2_{x,v}}^2\lesssim \|f^\eps_0\|_{L^2_{x,v}}^2+\sup_{t\in[0,T]}|M^\eps(t)|.
\end{equation}
On the other hand, similarly to \refe{2varM2}, we have
$$
\langle M^\eps,M^\eps\rangle(t)=\int_0^t (M|\varphi_1|^2-2\varphi_1M\varphi_1)(f^\eps(s),m^\eps(s))ds.
$$
(Note that there is no remaining terms here since $\varphi_2\equiv 0$, {\it cf.} the proof of \refe{2varM2} in Corollary~\ref{c9}.) In particular, by \refe{quadvarphi1l2} and the similar estimate for $M\varphi_1$, we have 
$$
|\langle M^\eps,M^\eps\rangle(t)|\lesssim \int_0^t \|f^\eps(s)\|_{L^2_{x,v}}^4ds.
$$
Since $M^\eps$ is a martingale with $\E M^\eps(t)=0$, Burkholder-Davis-Gundy inequality gives
\begin{equation}\label{afterDoob}
\E[\sup_{t\in[0,T]}|M^\eps(t)|^p|]\leq C_p \E |\langle M^\eps,M^\eps\rangle(T)|^{p/2}\lesssim \int_0^T\E\|f^\eps(s)\|_{L^2_{x,v}}^{2p}ds.
\end{equation}
By \refe{boundl2Meps}, $\E\|f^\eps(t)\|_{L^2_{x,v}}^{2p}\lesssim \E\|f^\eps_0\|_{L^2_{x,v}}^{2p}+\E[\sup_{t\in[0,T]}|M^\eps(t)|^p]$. Hence, by \refe{afterDoob},
$$
\E\|f^\eps(T)\|_{L^2_{x,v}}^{2p}\lesssim \E\|f^\eps_0\|_{L^2_{x,v}}^{2p}+ \int_0^T\E\|f^\eps(s)\|_{L^2_{x,v}}^{2p}ds.
$$
By Gronwall Lemma, we obtain $\E\|f^\eps(T)\|_{L^2_{x,v}}^{2p}\lesssim \E\|f^\eps_0\|_{L^2_{x,v}}^{2p}$. This actually holds true for any $t\in[0,T]$. Thus, using \refe{afterDoob} and then \refe{boundl2Meps} gives finally $\E[\sup_{t\in[0,T]}\|f^\eps(t)\|_{L^2_{x,v}}^{2p}]\lesssim \E\|f^\eps_0\|_{L^2_{x,v}}^{2p}$. \qed

\subsection{Tightness}\label{sec:tight}

\begin{proposition}[Tightness] Let $T>0$, $\eta>0$. Assume \refe{mmoments}, \refe{regNoise}, \refe{VanishFirstMoment}, \refe{PositiveSecondMoment} and assume that $(f^\eps_0)$ is bounded in $L^2$. Then $(\rho^\eps)$ is tight in $C([0,T];H^{-\eta}(\T^d))$.
\label{prop:tight}\end{proposition}

{\bf Proof:} Let $\varphi(\rho)=\rho$, or, more precisely (since the perturbed test-function method has been developed for real-valued, regularizing functions), define the test-function $\varphi_j$ as follows. Let $\{p_j;j\geq 1\}$ be a complete orthonormal system in $L^2(\T^d)$, let $\gamma> \max\{3,d\}$ and let 
$$
J=(\mathrm{Id}-\Delta_x)^{-1/2},
$$ 
where $\mathrm{Id}$ is the identity on $L^2(\T^d)$.  Note that $\|\rho\|_{H^{-\gamma}(\T^d)}=\|J^\gamma\rho\|_{L^2}$ and that $J^\gamma$ is Hilbert-Schmidt on $L^2(\T^d)$ since $\gamma>d$. We set $\varphi_j(\rho)=(J^\gamma\rho,p_j)$. It is clear that $\varphi_j$ is subquadratic (it is linear) and regularizing as in \refe{regularizingsubquadratic} (the operator $\nabla^3J^\gamma$ is of order $\leq 0$). Let $\varphi^\eps_j$ be the correction of $\varphi_j$ given by Corollary~\ref{c9}. Let 
$$
M^\eps_j(t)=\varphi^\eps_j(f^\eps(t),m^\eps(t))-\varphi^\eps_j(f^\eps_0,m^\eps(0))-\int_0^t\L^\eps\varphi^\eps_j(f^\eps(s),m^\eps(s))ds
$$
and
$$
\theta^\eps_j(t)=\varphi_j(\rho^0)+\int_0^t\L^\eps\varphi^\eps_j(f^\eps(s),m^\eps(s))ds+M^\eps_j(t).
$$
Then 
$$
\varphi_j(\rho^\eps(t))-\theta^\eps_j(t)=[\varphi_j(\rho^\eps(t))-\varphi_j^\eps(f^\eps(t),m^\eps(t))]-[[\varphi_j(\rho^\eps_0)-\varphi_j^\eps(f^\eps_0,m^\eps(0))]].
$$
By the estimates~\refe{corrsubq} on the correctors of $\varphi_j$ and the $L^2$-bounds of Proposition~\ref{prop:L2bound}, we deduce that
\begin{equation}\label{eqmodifid}
\E[\sup_{t\in[0,T]}|\varphi_j(\rho^\eps(t))-\theta^\eps_j(t)|]\lesssim\eps,
\end{equation}
where we write $a(j,\eps)\lesssim b(j,\eps)$ if there exists a constant $C$ depending on $\|a\|_{L^\infty(V)}$, on $T$, on the constant $C_*$ in \refe{regNoise} and on $\sup_{\eps>0}\|f^\eps_0\|_{L^2}$, but not on $\eps$ and $j$ such that $a(j,\eps)\leq C b(j,\eps)$.
Note also that, by \refe{corrsubq}, $\E\sup_{t\in[0,T]}|\theta^\eps_j(t)|\lesssim 1$, hence 
$$
\theta^\eps(t):=J^{-\gamma}\sum_{j\geq 1}\theta^\eps_j(t) J^\gamma p_j
$$
is a.s. well defined for all $t\in[0,T]$ in $H^{-\gamma}(\T^d)$ since the sum is convergent in $L^2(\T^d)$. By \refe{eqmodifid}, we obtain
\begin{equation}\label{estimmodifid}
\E[\sup_{t\in[0,T]}\|\rho^\eps(t)-\theta^\eps(t)\|_{H^{-\gamma}(\T^d)}]\lesssim\eps.
\end{equation}
Let $\eta>0$. Let 
$$
w(\rho,\delta):=\sup_{|t-s|<\delta}\|\rho(t)-\rho(s)\|_{H^{-\eta}(\T^d)}
$$
denote the modulus of continuity of a function $\rho\in C([0,T];H^{-\eta}(\T^d))$. Since the injection $L^2(\T^d)\subset H^{-\eta}(\T^d)$ is compact, the set
$$
K_{R}=\left\{\rho\in C([0,T];H^{-\eta}(\T^d)); \sup_{t\in[0,T]}\|\rho(t)\|_{L^2}\leq R; w(\rho,\delta)\leq \eps(\delta)\right\},
$$
where $R>$ and $\eps(\delta)\to 0$ when $\delta\to 0$, is compact in $C([0,T];H^{-\eta}(\T^d))$ (Ascoli's Theorem). By Prokhorov's Theorem, the tightness of $(\rho^\eps)$ will follow if we prove that, for all $\alpha>0$, there exists $R>0$, such that
\begin{equation}\label{tcompact}
\P(\sup_{t\in[0,T]}\|\rho^\eps(t)\|_{L^2}>R)<\alpha,
\end{equation}
and
\begin{equation}\label{cpctmodulus}
\lim_{\delta\to 0}\limsup_{\eps\to 0}\P(w(\rho^\eps,\delta)>\alpha)=0.
\end{equation}
The estimate \refe{tcompact} follows from the $L^2$-bound of Proposition~\ref{prop:L2bound} by the estimate
$$
\P(\sup_{t\in[0,T]}\|\rho^\eps(t)\|_{L^2}>R)\leq\frac{1}{R}\E[\sup_{t\in[0,T]}\|\rho^\eps(t)\|_{L^2}].
$$
Similarly, we will deduce \refe{cpctmodulus} from a bound on $\E w(\rho^\eps,\delta)$ for $\delta>0$. Actually, by the $L^2$-bound of Proposition~\ref{prop:L2bound} and by interpolation, we have
$$
\E\sup_{|t-s|<\delta}\|\rho(t)-\rho(s)\|_{H^{-\eta^\flat}(\T^d)}\leq \E\sup_{|t-s|<\delta}\|\rho(t)-\rho(s)\|^\sigma_{H^{-\eta^\sharp}(\T^d)}
$$
for a certain $\sigma>0$ if $\eta^\sharp>\eta^\flat>0$. Therefore it is indeed sufficient to work with $\eta=\gamma$. Besides, by \refe{estimmodifid}, it is sufficient to obtain an estimate on the increments of $\theta^\eps$. By definition
$$
\theta^\eps_j(t)-\theta^\eps_j(s)=\int_s^t\L^\eps\varphi^\eps_j(f^\eps(\sigma),m^\eps(\sigma))d\sigma+M^\eps_j(t)-M^\eps_j(s),
$$
for $0\leq s\leq t\leq T$. By the $L^2$-bound of Proposition~\ref{prop:L2bound} and by \refe{corrsubq}, we have
$$
\E\left|\int_s^t\L^\eps\varphi^\eps_j(f^\eps(\sigma),m^\eps(\sigma))d\sigma\right|^4\lesssim |t-s|^4.
$$
By Burkholder-Davis-Gundy inequality,
$$
\E|M^\eps_j(t)-M^\eps_j(s)|^4\lesssim \E|\<M^\eps_j,M^\eps_j\>(t)-\<M^\eps_j,M^\eps_j\>(s)|^2,
$$
where $\<M^\eps_j,M^\eps_j\>$ is the quadratic variation of $M^\eps_j$. By \refe{2varM2} and the $L^2$-bound of Proposition~\ref{prop:L2bound}, we obtain
$$
\E|M^\eps_j(t)-M^\eps_j(s)|^4\lesssim |t-s|^2.
$$
Finally, we have $\E|\theta^\eps_j(t)-\theta^\eps_j(s)|^4\lesssim |t-s|^2$, and thus 
$$
\E\|\theta^\eps(t)-\theta^\eps(s)\|_{H^{-\gamma}(\T^d)}^4\lesssim |t-s|^2.
$$
It follows (by the Kolmogorov's criterion) that, for $\alpha<1/2$, 
$$
\E\|\theta^\eps\|^4_{W^{\alpha,4}(0,T;H^{-\gamma}(\T^d))}\lesssim 1.
$$ 
By the embedding
$$
W^{\alpha,4}(0,T;H^{-\gamma}(\T^d))\subset C^{0,\mu}([0,T];H^{-\gamma}(\T^d)),\quad \mu<\alpha-\frac{1}{4},
$$
we obtain $\E w(\theta^\eps,\delta)\lesssim \delta^\mu$ for a certain positive $\mu$. This concludes the proof of the proposition. \qed

\subsection{Convergence}\label{sec:convergence}

We conclude here the proof of Theorem~\ref{th:DiffLimStochastic}. Fix $\eta>0$. By Proposition~\ref{prop:tight}, there is a subsequence still denoted by 
$(\rho^\eps)$ and a probability measure $P$ on $C([0,T];H^{-\eta}(\T^d))$ such that 
$$
P^\eps\to P\quad\mbox{weakly on}\quad C([0,T];H^{-\eta}(\T^d))
$$
where $P^\eps$ is the law of $\rho^\eps$. We then show that $P$ is a solution of the martingale problem, with a set of test
functions specified below, associated to  the limit equation \eqref{LimSPDE}.
\medskip

By Skohorod representation Theorem \cite{Billingsley99}, and since $C([0,T];H^{-\eta}(\T^d))$ is 
separable, there exists a new probability space $(\widetilde{\Omega},\widetilde{\mathcal{F}},\widetilde{\P})$ and some random va\-ria\-bles 
\begin{equation*}
\widetilde{\rho}^\eps,\widetilde{\rho}\colon\widetilde{\Omega}\to C([0,T];H^{-\eta}(\T^d)),
\end{equation*}
with respective law $P^\eps$ and $P$ such that $\widetilde{\rho}^\eps\to\widetilde{\rho}$ in $C([0,T];H^{-\eta}(\T^d))$, $\widetilde{\P}$ a.s. 
\medskip

Let $\varphi\in C^3(L^2(\T^d))$ be regularizing and subquadratic according to \refe{regularizingsubquadratic}. By Corollary \ref{c9} and the $L^2$-bound of Proposition~\ref{prop:L2bound}, we have  for $0\le s_1\le \dots \le s_n\le s\le t$ and
$\psi\in C_b((L^2_{x,v})^n)$,
\begin{equation}\label{e41}
\ds \left|\E\left[\left(\varphi(\rho^\eps(t))-\varphi(\rho^\eps(s))-\int_s^t\L\varphi(\rho^\eps(\sigma))d\sigma\right)\psi(\rho^\eps(s_1),\dots,\rho^\eps(s_n))\right]\right|\le C \eps
\end{equation}
with a constant $C$ depending on the constant $C_*$ in \refe{regNoise}, on $C_\varphi$, on $\sup_{\eps>0}\|f^\eps_0\|_{L^2}$ and on 
the supremum of $\psi$. Since $\rho^\eps$ and $\widetilde\rho^\eps$ have the same law, this is still true if $\rho^\eps$
is replaced by $\widetilde \rho^\eps$. 
Assume furthermore that $\varphi$ is bounded and continuous from $H^{-\eta}(\T^d)$ into $\R$, then it is easy to take the limit $\eps\to 0$ in \eqref{e41}
 and to obtain
\begin{equation}
\label{e42}
\widetilde{\E}\left\{\left[\varphi(\widetilde{\rho}(t))-\varphi(\widetilde{\rho}(s))-\int_s^t\L\varphi(\widetilde{\rho}(\sigma))d\sigma\right]\psi(\widetilde{\rho}(s_1),\dots,\widetilde{\rho}(s_n))\right\}=0.
\end{equation}
The additional hypothesis on $\varphi$ can be relaxed. Indeed, thanks to Proposition \ref{prop:L2bound}, we can approximate every subquadratic and regularizing functions by functions in $C_b(H^{-\eta}(\T^d))$ which are subquadratic and regularizing with a uniform constant in \refe{regularizingsubquadratic} and which converge pointwise.
\medskip

We have thus proved that $P$ solves the martingale problem associated to $\L$ with 
subquadratic and regularizing test functions. In particular, for all such 
$\varphi$:  
\begin{equation}
M_\varphi(t)= \varphi(\rho(t))-\int_0^t\L\varphi(\rho(s))ds, \quad t\ge 0,
\label{tildemartingale}\end{equation}
is a martingale
with respect to the filtration ${\mathcal{F}}_s$ generated by $(\rho(s))$. The quadratic variation of $M_\varphi$ is ({\it cf} \refe{2varM})
$$
\<M_\varphi,M_\varphi\>(t)=\L|\varphi|^2-2\varphi\L\varphi.
$$
Let $D\varphi(\rho)\otimes D\varphi(\rho)$ denote the bilinear form 
$$
(h,k)\mapsto (h,D\varphi(\rho))(k,D\varphi(\rho))
$$ 
on $L^2(\T^d)$. By \refe{devLQ}, we have
\begin{align*}
\L\varphi^2(\rho)-2\varphi(\rho(s))\L\varphi(\rho(s))&=\frac12\mathrm{Trace}[\rho Q^{1/2})D\varphi(\rho)\otimes D\varphi(\rho)\rho Q^{1/2})^*]\\
&=\left\| \rho Q^{1/2} D\varphi(\rho)\right\|_{L^2}^2.
\end{align*}
We deduce that 
$$
M(t)= \rho(t)-\rho(0)-\int_0^t {\mathcal A}\rho(s)+\frac12F\rho(s)ds, \; t\ge 0,
$$ 
is a martingale with quadratic variation $\ds\int_0^t \rho(s)Q^{1/2}\left(\rho(s)Q^{1/2}\right)^*ds$. Thanks to martingale representation results (see for instance \cite{DaPratoZabczyk92}), up to
a change of probability space, there exist a cylindrical Wiener process $W$ such that 
$$
\rho(t)-\rho(0)-\int_0^t {\mathcal A}\rho(s)+\frac12F\rho(s)ds=\int_0^t \rho(s)Q^{1/2}dW(s),\quad t\ge 0.
$$
It is well known that this equation has a unique solution with paths in the space $C([0,T];H^{-\eta}(\R^d))$. 
This can be shown for instance by energy estimates using It\^o formula after a suitable
regularization argument. 
Moreover pathwise uniqueness implies 
uniqueness in law and we deduce that $P$ is the law of this solution and is uniquely determined. 
Finally, by uniqueness of the limit, the whole sequence $(P^\eps)$ converges to $P$ weakly in the space of 
probability measures on $C([0,T];H^{-\eta}(\R^d))$. 

\def\cprime{$'$} \def\cprime{$'$}
\providecommand{\bysame}{\leavevmode\hbox to3em{\hrulefill}\thinspace}
\providecommand{\MR}{\relax\ifhmode\unskip\space\fi MR }
\providecommand{\MRhref}[2]{%
  \href{http://www.ams.org/mathscinet-getitem?mr=#1}{#2}
}
\providecommand{\href}[2]{#2}

\end{document}